\documentclass[12pt]{article}
\usepackage{amsfonts, amsmath, amssymb, amsthm}
\usepackage{amscd, eucal}

\input xy
\xyoption{all}

\usepackage{layout}
\usepackage{graphicx}

\newtheorem{theorem}{Theorem}
\newtheorem{propos}[theorem]{Proposition}

\newtheorem{rem}[theorem]{Remark}
\newtheorem{defi}[theorem]{Definition}
\newtheorem{qes}[theorem]{Question}
\newtheorem{lemma}[theorem]{Lemma}
\newtheorem{example}[theorem]{Example}

\newcommand{\chapter}[1]{{\bf\Large #1}}

\def\beq{\begin{equation}}
\def\eeq{\end{equation}}
\def\blm{\begin{lemma}}
\def\elm{\end{lemma}}
\def\bdf{\begin{defi}}
\def\edf{\end{defi}}
\def\btm{\begin{theorem}}
\def\etm{\end{theorem}}
\def\bpp{\begin{propos}}
\def\epp{\end{propos}}
\def\bQ{\begin{qes}}
\def\eQ{\end{qes}}
\def\ben{\begin{enumerate}}
\def\een{\end{enumerate}}
\def\into{\hookrightarrow}
\def\bex{\begin{example}}
\def\eex{\end{example}}
\def\brm{\begin{rem}}
\def\erm{\end{rem}}

\def\lb{\label}
\def\ger{\mathfrak}

\input{epsf.sty}

\DeclareMathAlphabet{\mathscr}{OT1}{pzc}%
                                 {m}{it}

\newcommand{\Z}{{\mathbb Z}}

\newcommand{\Cc}{{\mathbb C}}

\newcommand{\Q}{{\mathbb Q}}

\begin{document}

\title{Nilpotent slices and Hilbert schemes}

\author{Craig Jackson}
\date{}

\maketitle

\pagestyle{plain}



\begin{abstract}
We construct embeddings 
$\mathcal{Y}_{n,\tau} \to \mbox{Hilb}^n (\Sigma_{\tau})$
for each of the classical Lie algebras  $\ger{sp}_{2m}(\Cc)$, 
$\ger{so}_{2m}(\Cc)$, and $\ger{so}_{2m+1}(\Cc)$. 
The space $\mathcal{Y}_{n,\tau}$ is the fiber over a point
$\tau \in \ger h / W$ of the restriction of the adjoint
quotient map $\chi : \ger g \to \ger h /W$ to a suitably chosen 
transverse slice of a nilpotent orbit. These embeddings
were discovered for $\ger{sl}_{2m}(\Cc)$ by Ciprian Manolescu.
They are related to the symplectic link homology of Seidel and Smith.
\end{abstract}

\section{Introduction}\lb{intro}


Khovanov cohomology is a link invariant that takes the form of a 
bi-graded Abelian group $Kh^{i,j}(L)$ (\cite{Kho}). It is a categorification of
the Jones polynomial in the sense that, up to normalization and
change of variables, the Jones polynomial is given as the graded Euler
charateristic of $Kh^{i,j}(L)$:
\[
V_L (t)= (q+q^{-1})^{-1}
\sum_{i,j} (-1)^{i} q^j \dim(Kh^{i,j} \otimes \Q)|_{q = -t^{1/2}}.
\]
Khovanov cohomology is
known to be a strictly stronger invariant than the Jones polynomial and is,
by definition, able to be computed combinatorially by means of long 
exact sequences. However, unlike some other knot invariants 
(the Alexander polynomial, for instance) 
the geometric meaning of both the Jones polynomial and 
Khovanov cohomology has not been fully understood.

Recently, Seidel and Smith \cite{SS} proposed a geometric
interpretation of Khovanov cohomology in terms of symplectic geometry and 
Lagrangian Floer cohomology.
They start by presenting a link $L$ as the closure of an
$m$-strand braid $\beta \in B_m$. Appending $m$ trivial stands
gives $\beta \times 1^m \in B_{2m}$ which can be viewed as
a path in the configuration space 
$\mbox{Conf}_{2m}(\Cc) = (\Cc^{2m} - \Delta)/S_{2m}$.
They construct a symplectic fibration over $\mbox{Conf}_{2m}(\Cc)$
as a slice $\mathcal{S}_m$ 
of a 
distinguished nilpotent orbit in $\ger{sl}(2m)$. The fiber over
each $\tau \in \mbox{Conf}_{2m}(\Cc)$ is a symplectic manifold
$\mathcal{Y}_{m,\tau}$ and the monodromy along $\beta$ gives rise
to a symplectic automorphism $\phi_{\beta}$ of $\mathcal{Y}_{m,\tau}$. 
They chose a distinguished
Lagrangian submanifold $\mathcal L \subset \mathcal{Y}_{M,\mu}$
and apply the monodromy map to obtain another Lagrangian
$\phi_{\beta} \mathcal L$. To this geometric data they apply
the theory of Lagrangian Floer cohomology and define the
\emph{symplectic Khovanov cohomology} to be
\[
Kh_{\mbox{\scriptsize{symp}}}^{k} (L) = 
      HF^{k+m+\omega} (\mathcal L, \phi_{\beta} \mathcal L)
\] 
where $\omega$ is the writhe of the link $L$ (the signed number of
crossings).

They show that, up to an isomorphism of graded abelian groups,
$Kh_{\mbox{\scriptsize{symp}}}$ is invariant under the
markov moves and so defines a link invariant. Furthermore,
they conjecture that symplectic Khovanov cohomology is the same as 
the normal bi-graded cohomology after a collapsing of the bi-grading:
\[
Kh_{\mbox{\scriptsize{symp}}}^{k} (L) = \bigoplus_{i+j=k} Kh^{i,j}(L).
\]

Ciprian Manolescu \cite{Man} has discovered an interesting property
of the symplectic manifolds $\mathcal{Y}_{m,\tau}$. Namely, that
they can be embedded 
as an open dense subset of the Hilbert scheme $\mbox{Hilb}^m (S_{\tau})$
of $m$ points on the complex surface 
$S_{\tau} = \{(u,v,z) \in \Cc^3 \mid u^2 + v^2 + P_{\tau}(z) = 0\}$.
Here $P_{\tau}(t)$ is the unique monic polynomial with roots
given by $\tau$ and $S_{\tau}$ is a deformation of 
the $A_{2m-1}$ Kleinian singularity.

This embedding is interesting for a number of reasons. First, it gives 
a more concrete form to the geometrical construction of Seidel and Smith.
Composing with the Hilbert-Chow morphism
$\pi_{HC} : \mbox{Hilb}^m (S_{\tau}) \to \mbox{Sym}_m (S_{\tau})$ gives 
nice holomorphic coordinates on $\mathcal{Y}_{m,\tau}$. These coordinates,
as well as a natural $\Cc^*$-action on the nilpotent slice, 
can be used to define a Kahler metric 
on $\mathcal{S}_m$.
This metric descends to give a symplectic structure on each
$\mathcal{Y}_{m,\tau}$ which is shown to be a deformation of
the symplectic structure given by Seidel and Smith. 

Secondly, in 
Manolescu's setup, the Hilbert-Chow
morphism 
can be used to concretely describe the Lagrangians $\mathcal L$ and
$\phi_{\beta} \mathcal{L}$
as a product of $m$ 2-spheres living in $S_{\tau}$. Because
of this, the intersection 
$\mathcal L \cap \phi_{\beta} \mathcal{L}$ can be explicitly described. 
This in turn gives
a set of generators of the Seidel-Smith cohomology. Manolescu is 
then able to set up a natural correspondence between these generators
and the set of intersection points in Bigelow's picture of the Jones
polynomial \cite{Big}.

Lastly, the embedding $\mathcal{Y}_{m,\tau} \into \mbox{Hilb}^m (S_{\tau})$
is interesting for its own sake. Both spaces are examples of
quiver varieties. Quiver varieties are built by attaching pairs of
Hermetian vector spaces to each node of a finite oriented 
graph. In fact, the varieties 
$\mathcal{Y}_{m,\tau}$ and  $\mbox{Hilb}^m (S_{\tau})$
are constructed from the same quiver, but with different data on the
nodes. The embedding $\mathcal{Y}_{m,\tau} \into \mbox{Hilb}^m (S_{\tau})$
is then an open holomorphic embedding from one quiver variety
into another.

In this paper we construct Manolescu-type embeddings for the remaining
classical Lie algebras: $\ger{sp}_{2m}(\Cc)$, $\ger{so}_{2m}(\Cc)$, 
and $\ger{so}_{2m+1}(\Cc)$.  Namely, we show that for certain 
``degree two'' nilpotent orbits, a slice can be chosen  such that
$(i)$ the restiction of the adjoint 
quotient to the slice is a smooth fiber bundle and
$(ii)$ each fiber can be embedded holomorphically as a dense open subset 
of the Hilbert scheme of $m$ points on a complex
surface $\Sigma$. In all cases, the slice in question needs to be chosen
carfully so that the open embedding into the Hilbert scheme can be achieved.
That is, these slices are not
the standard Slodowy slices given by the Jacobson-Morozov lemma, even though
they are similar in a few important ways.
Also, whereas in the type-$A$ case the complex surface $\Sigma$
is a deformation of a type-$A$ singularity,
in the case of the other complex semisimple algebras, the
surface is a deformation of a type-$D$ singularity.

\section{The Geometry of Nilpotent Orbits}\lb{geom}

The embedding in the type-$A$ case is constructed as an algebraic morphism
by first specifying certain ``second order'' nilpotent orbits 
and then finding slices of these orbits 
whose elements have a nice form to their characteristic
polynomials. The same basic
recipe holds in the other cases as well, though there are some
interesting diferences. Note, however, that in all cases the transverse
slices have to be chosen carefully, as the standard Slodowy slices
given by the Jocobson-Morozov lemma have rather complicated characteristic
polynomials. 
Before we get to the specifics,
however, we sketch out some of the results we will need regarding the
geometry of nilpotent orbits and arbitrary transverse slices. So as to 
keep things as concise as possible we restrict our attention 
to the Lie algebras of types $B$, $C$, and $D$.

\subsection{The Adjoint Quotient}

Let $\mathfrak{g}$ be a complex semisimple Lie algebra
of type $B$, $C$, or $D$. That is, 
$\ger g$ is $\ger{sp}_{2m}(\Cc)$, 
$\ger{so}_{2m}(\Cc)$, or
$\ger{so}_{2m+1}(\Cc)$. Let $G$ be the appropriate 
complex Lie group acting on $\ger g$
by the adjoint action. For any $Y$ in 
$\ger g$, denote the adjoint orbit of $Y$ by 
 $\mathcal{O}_Y =  Ad(G) Y$.
Let $\ger h$ be the standard
Cartan subalgebra of all diagonal matrices in $\ger g$, 
and let $W = N_G (\ger h ) / C_G (\ger h )$ 
be the the Weyl group. In each case, the 
Cartan $\ger h$ is isomorphic to the space
$D_m = \{ \mbox{diag}(x_1, \dots, x_m) \mid x_i \in \Cc \}$ of
all $m \times m$ complex diagonal matrices since, physically, the space 
$\ger h$ is the ``skew-diagonal'' inside $D_m \oplus D_m$. 
That is, $\ger h = \{ (x, -x) \mid x \in D_m \}$ for types $C$ and $D$
and  $\ger h = \{ (0, x, -x) \mid x \in D_m \}$
for type-$B$.

For the algebras $\ger{sp}_{2m}(\Cc)$ and $\ger{so}_{2m+1}(\Cc)$ the
Weyl group $W$ is the same. It 
is the semi-direct product of 
the symmetric group $S_m$ with $(\Z /2)^m$, where
$S_m$ acts on $\ger h \cong D_m$ 
by permuting the coordinates,
and the $i$th factor of $(\Z / 2)^m$ acts by -1 on the $i$th coordinate.
We define a map 
\beq
\epsilon : \ger h / W \to \mbox{Sym}_m (\Cc)
\eeq
by $\epsilon : [ \mbox{diag}(x_1, \dots, x_m) ] 
\mapsto \{x_1^2, \dots, x_m^2 \}$. 
By $\mbox{Sym}_m (\Cc)$ we mean  the \emph{symmetric product} 
of $m$ copies of $\Cc$. That is, $\mbox{Sym}_m (\Cc) = \Cc^m / S_m$.
For type $B$ and $C$ algebras the map $\epsilon$ is an isomorphism.

For the algebra $\ger{so}_{2m}(\Cc)$ the Weyl group is
isomorphic to the semi-direct product of
the symmetric group $S_m$ with $(\Z /2)^{\binom{m}{2}}$, where
$S_m$ acts on $\ger h \cong D_m$
by permuting the coordinates,
and the $(i,j)$th factor of $(\Z / 2)^{\binom{m}{2}}$ acts by 
-1 on both the $i$th and $j$th coordinate.
Suppose, in this case, we define a $\Z / 2$-action on
$\ger h / W$ that multiplies a single coordinate by -1. 
The set of fixed points of this action corresponds to the set
of all elements of $\ger h$ that have at least one zero
eigenvalue. Because of this, 
let us write the set of $\Z / 2$-fixed points as
$(\ger h / W)^0$.
Hence, in the type-$D$ case, if we take the map
$\epsilon : \ger h / W \to \mbox{Sym}_m (\Cc)$ 
given above, then we see that this map has fibers
given by the $\Z / 2$ action. Thus, this map is an isomorphism when
restricted to $(\ger h / W)^0$. Also, it is clear that 
the image of $(\ger h / W)^0$
under $\epsilon$ is the set $\mbox{Sym}_{m}^0 (\Cc)$: the set of all
elements of the symmetric product having at least one point equal 
to zero.

Now, for any $Y \in \ger g$ let us write $Y_s$ for its semisimple
part. The space $\ger h / W$ is in natural bijective correspondence with
set of all semisimple $G$-orbits 
by the map 
$\mu : [H] \to \mathcal{O}_H$ (\cite{CM} 2.2). 
We define the \emph{adjoint quotient map} 
$\chi : \ger g \to \ger h / W$ by
$\chi : Y \mapsto {\mu}^{-1}(\mathcal{O}_{Y_s})$.
Now, $Y$ and $Y_s$ have the same set of eigenvalues so by composing $\chi$
with the map $\epsilon :\ger h / W \to \mbox{Sym}_m (\Cc)$ we
can interpret $\chi$ as the map that takes $Y$ to its
set of generalized eigenvalues.

Let us define
$\bar \chi : \ger g \to \mbox{Sym}_{2m} (\Cc)$ by setting $\bar \chi (Y)$ equal
to the unordered $2m$-tuple of roots of $\chi_{Y} (t)$, the characteristic
polynomial of $Y$. (In the type-$B$ case we throw away the trivial 
$t=0$ that is a root of every characteristic polynomial.)

We may also define an embedding $\mbox{Sym}_m(\Cc) \into \mbox{Sym}_{2m}(\Cc)$
by $\{ \mu_1, \dots, \mu_m \} 
 \mapsto \{ \lambda_1 , -\lambda_1 , \dots , \lambda_m , -\lambda_m \}$, where
$\lambda_i$ is any square root of $\mu_i$. 
We have already seen that the eigenvalues of any element $Y \in \ger g$
occur as opposite pairs $(x, -x)$ with the addition of one trivial root
in the type-$B$ case.  This means that 
the image of this embedding  $\mbox{Sym}_m(\Cc) \into \mbox{Sym}_{2m}(\Cc)$
is identical to the image of $\bar \chi$.

All of this fits together into the following commutative diagram:
\[
\xymatrix{ \ger g \ar[rr]^(.4){\bar \chi} \ar[d]_{\chi} 
                     && \mbox{Sym}_{2m} (\Cc) \ar@{<-_{)}}[d] \\
           \ger h /W \ar[rr]^(.4){\epsilon} && \mbox{Sym}_m (\Cc)      }.
\]

\subsection{Transverse Slices}

A \emph{transverse slice} to the orbit $\mathcal{O}_X$
is a complex submanifold $\mathcal S \subseteq \ger g$
whose tangent space at $X$ is
complementary to $T_X ( \mathcal{O}_X )$. Now the tangent
space to the orbit is given by
\[
T_X ( \mathcal{O}_X ) = [ \ger g , X].
\] 
So if the transverse slice 
$\mathcal S$ happens to be an affine subspace of $\ger g$,
with $X \in \mathcal S$, then $\mathcal S$ will be of the
form 
\[
\mathcal S = X + V
\]
for some vector subspace $V \subseteq \ger g$
that is complementary to $[ \ger g , X]$. 

We have the following local result on 
transverse slices (\cite{Slo1} 5.1, 
\cite{SS} 2A):

\blm\lb{lemma5}
(i) For all $Y \in \mathcal S$ sufficiently close to $X$, the
intersection $\mathcal S \cap \mathcal{O}_Y$ is transverse at $Y$.
(ii) For all $Y \in \mathcal S$ sufficiently close to $X$, 
$Y$ is a critical point of $\chi | {\mathcal S}$ if and only if
it is a critical point of $\chi $.
(iii) Any two transverse slices at $X$ are locally isomorphic
by an orbit preserving isomorphism.
\elm

Define a vector field $\xi$ on $\ger g$ by 
\[
\xi_Y = Y - \frac{1}{2}[H,Y].
\] 
This is the infinitesimal generator of a linear $\Cc^*$-action $\lambda$
on $\ger g$. It is easy to check that
this action is given by
\beq\lb{cstar}
\lambda_r (Y) = r \mbox{Ad}(r^{H/2})Y 
           = r \exp (-\log(r) H/2) Y \exp(\log(r) H/2).
\eeq
Notice that $\lambda$ preserves orbits in the sense
that $\lambda_r (\mathcal{O}_Y) = \mathcal{O}_{\lambda_r (Y)}$. 
On the level of
characteristic polynomials $\lambda_r$ multiplies every
eigenvalue by $r$. Also, since $\xi_X = 0$, we see
that $X$ is a fixed point of $\lambda$.

\bdf
A slice $\mathcal S$ transverse to the orbit $\mathcal{O}_X$
at the point $X$ is said to be $\lambda$-invariant if it is
invariant under the $\Cc^*$-action given by \eqref{cstar}.
\edf

\begin{example}
If $X \in \ger g$ is a nilpotent, then there is a canonical way 
of constructing transverse slices using the Jacobson-Morozov Lemma. This
lemma says that there are elements $H, N^+, N^- \in \ger g$, 
with $N^+ = X$, such that
\begin{align}
[H,N^+] &= 2N^+ & [H,N^-] &= -2N^- & [N^+,N^-] &= H.
\end{align}

This produces a splitting $\ger g = [ \ger g , N^+ ] \oplus {\ger g}_{N^-}$,
which shows that the space 
$\mathcal{S}^{\mbox{\scriptsize{JM}}} = X + {\ger g}_{N^-}$ 
is a transverse slice to $\mathcal{O}_X$. The space 
$\mathcal{S}^{\mbox{\scriptsize{JM}}}$ is called the
\emph{JM-slice} or the \emph{Slodowy slice}
(\cite{Slo1} 7.4). There are many choices
of JM triples for a fixed $X$, but they are all conjugate
to eachother by an element of the stabilizer $G_X$.

Now if $Y$ is an element of ${\ger g}_{N^+}$, then so is $[H,Y]$. Thus,
$\xi$ restricts to a vector field on $\mathcal{S}^{\mbox{\scriptsize{JM}}}$,
which means that $\mathcal{S}^{\mbox{\scriptsize{JM}}}$ is invariant
under $\lambda$.
\end{example}

The Slowdowy slice is the prototypical example of 
a $\lambda$-invariant slice. These are not the only examples, however,
as we will see later. 
For all $\lambda$-invariant slices we have the following result,
which is an improvement over Lemma \ref{lemma5} (cf. Lemmas 15 and 16 in 
\cite{SS}).

\blm\lb{lemma15}
Let $\mathcal S$ be a $\lambda$-invariant slice at $X$. Then
(i) the intersection of $\mathcal S$ with any adjoint orbit is transverse.
(ii) A point of $\mathcal S$ is a critical point of $\chi | \mathcal S$
if and only if it is a critical point of $\chi$. 
(iii) Any two $\lambda$-invariant slices are isomorphic by a 
$\Cc^*$-equivariant, orbit preserving isomorphism.
\elm

\begin{proof}
The  $\Cc^*$-action shrinks $\mathcal S$ to $X$. That is, for any
$Y \in \mathcal S$ we have $\lim_{r \to 0} \lambda_r (Y) = X$. Hence, the
local properties of the slice given in Lemma \ref{lemma5} will 
extend to the entire slice. 
\end{proof}

\subsection{Simultaneous Resolution}

Consider the open subset $\mbox{Conf}_m^{\,*} (\Cc) = \mbox{Conf}_m (\Cc^*) 
        \subset \mbox{Sym}_m (\Cc) \cong \ger h /W$
consisting of all sets of $m$ distinct, nonzero points in $\Cc$. Also,
let us define $\mbox{Conf}_m^{\,0} (\Cc)$ to be the subset of
$\mbox{Conf}_m (\Cc)$ consisting of all sets of $m$ distinct points
in $\Cc$ with one point equal to 0. Then clearly we have
$\mbox{Conf}_m (\Cc) = \mbox{Conf}_m^{\,*} (\Cc) \cup \mbox{Conf}_m^{\,0} (\Cc)$
and $\mbox{Conf}_m^{\,0} (\Cc) \cong \mbox{Conf}_{m-1}^{\,*} (\Cc)$.   

\brm
Note that in \cite{Man} and \cite{SS}, the notation 
$\mbox{Conf}_m^{\,0} (\Cc)$ stands for the set of all 
unordered $m$-tupiles that sum to zero.
\erm

For $Y \in \ger g$, suppose $\epsilon( \chi (Y))$ belongs to
$\mbox{Conf}_m^{\,*} (\Cc)$, 
then the eigenvalues of $Y$ must all be distinct. Therefore, $Y$ must be 
semisimple and since $\chi$ is $G$-invariant we can assume
$Y \in \ger h$ so that $Y$ is a regular point of the projection
$\ger h \to \ger h / W$. This implies that $Y$ is also a regular point of
$\chi : \ger g \to \ger h / W$. 

In the type-$D$ case we have a larger class of regular points since
the points in $\mbox{Conf}_m^{\,0} (\Cc)$ also correspond to semisimple
elements of $\ger g$ whose image in 
$\ger h$ do not lie in any root hyperplane.

To simplfy our discussion somewhat we define 
a subest $R_{\ger g} \subseteq \ger h / W$ as follows:
In the type-$B$ and type-$C$ case we set
$R_{\ger{g}}$ equal to the preimage of 
$\mbox{Conf}_m^{\,*} (\Cc)$ under the isomorphism $\epsilon$. 
In the type-$D$ case we have instead
$R_{\ger{so}_{2m}(\Cc)} = \epsilon^{-1}(\mbox{Conf}_m (\Cc))$.
Thus, for any $\tau \in R_{\ger g}$ we have that $\chi^{-1}(\tau)$
is a smooth submanifold of $\ger g$.

We  use Grothendieck's simultaneous resolution of
the adjoint quotient to improve on this result.

\bpp\lb{bundle1}
The adjoint quotient map $\chi : \ger g \to \ger h / W$ is
a smooth fiber bundle when restricted to $R_{\ger g}$.
\epp
\begin{proof}
Let $B$ be the standard Borel subgroup of $G$ and let 
$\ger b \subseteq \ger g$ be its
Lie algebra (the standard Borel subalgebra). 
Let $\ger{n} = [\ger b , \ger b ]$ be the nilradical. It is
the span of all root vectors
$X_{\alpha}$ for some positive root $\alpha$ and we have 
$\ger b = \ger h + \ger n$. Let $\mathcal B$ denote the space of all
Borel subalgabras of $\ger g$. Then $\mathcal B$ is a smooth
algebraic $G$-variety and is isomorphic to $G/B$ 
by $gB \mapsto Ad(g) \ger b$ (\cite{GC} 3.1). 

Let us denote
the projection of $\ger b$ onto $\ger h$ by $\pi_{\ger h}$.
Define $\tilde{\ger g}$ to be the fibered product
$G \times_B \ger b = (G \times \ger b)/B$, where the $B$ action is
defined by $b(g,Y) = (gb^{-1},Ad(b)Y)$. 
The space $\tilde{\ger g}$ is isomorphic to the incidence variety
$\{ (Y, \ger c) \in \ger g \times \mathcal B \mid Y \in \ger c \}$
by the $G$-equivariant map that takes $(g, Y)$ to $( Ad(g)Y, gB)$.

Define maps
$\tilde{\chi} : \tilde{\ger g} \to \ger h$ and 
$\psi : \tilde{\ger g} \to \ger g$ by
\begin{align}
\tilde{\chi} &: (g, Y) \mapsto \pi_{\ger h}(Y) \\
\psi &: (g,Y) \mapsto Ad(g)Y 
\end{align} 

The map $\tilde{\chi}$ is well defined since $B$ acts trivially on
$\ger h$ modulo $\ger n$ (\cite{Slo1} 4.3) and it factors through
the map $\tilde{\ger g} \to G/B \times \ger h$ given by
$(g, Y) \mapsto (gB, \pi_{\ger h}(Y))$. 
Thus, the mapping $\tilde{\chi} : \tilde{\ger g} \to \ger h$ 
has the natural structure
of a differentiable fiber bundle with fibers equal to
$G \times_{B} \ger n$.

The following commutative diagram 
is a simultaneous resolution of the adjoint quotient
map (cf. \cite{Slo1} 4.7, 
\cite{GC} 3.1)
\beq\lb{simres}
\xymatrix{\tilde{\ger g} \ar[rr]^(.4){\phi} \ar[d]_{\tilde{\chi}}
                     && \ger g \ar[d]_{\chi} \\
           \ger h  \ar[rr]^(.4){} && \ger h / W }.
\eeq

This means that for all $\mu \in \ger h$, the restriction
$\psi_{\mu} : \tilde{\chi}^{-1} (\mu) \to \chi^{-1}([\mu])$
is a resolution of singularities. In case $\chi^{-1}([\mu])$
is already smooth, the map $\psi_{\mu}$ is an isomorphism 
(cf. \cite{Slo1} 4.1).

For any $\tau \in R_{\ger g}$, we have already seen that $\chi^{-1}(\tau)$
is smooth. Also, the quotient map $\ger h \to \ger h / W$ is a smooth
$|W|$-sheeted covering when restricted to $R_{\ger g}$. This proves the 
proposition.
\end{proof}

We can restrict everything in diagram \eqref{simres}
to the transverse slice $\mathcal S$. What we 
obtain is a simultaneous resolution of $\chi | \mathcal{S}$ 
(\cite{Slo1} 5.2)
\beq\lb{simresslice}
\xymatrix{\tilde{\mathcal S} \ar[rr]^(.4){\phi | \tilde{\mathcal S}} 
                             \ar[d]_{\tilde{\chi}}
                     && \mathcal S \ar[d]_{\chi | \mathcal{S}} \\
           \ger h  \ar[rr]^(.4){} && \ger h / W }.
\eeq
If the slice happens to be 
$\lambda$-invariant, then from Lemma \ref{lemma15} we know that 
regular points of $\chi$ are
still regular after we restrict to $\mathcal S$. Hence,  
Proposition \ref{bundle1} carries over to give the following:

\bpp\lb{fiberbundle}
The restriction of the adjoint quotient map 
$\chi | \mathcal{S} : \mathcal S \to \ger h / W$ to 
a $\lambda$-invariant slice is
a smooth fiber bundle when restricted to $R_{\ger g}$.
\epp

We will generically denote the fiber of 
$\chi | \mathcal S$ over $\tau \in \ger h /W$
by $\mathcal{Y}_{\tau}$. We will see in the sections that follow
that it is an algebraic variety defined 
by $m = \mbox{rank } \ger g$ equations in the
coordinates of $\mathcal S$.

The last thing we will do in this section is mention a result
of Brieskorn and Slodowy (cf. \cite{Bri}, \cite{Slo1} 6.2, 8.7)  
relating nilpotent slices to simple singularities. 

\bpp\lb{Slodowy}
Let $X$ be a subregular element of a semsisimple complex Lie algebra
and let $\mathcal{S}$ be a JM-slice to the orbit of $X$. Then the
intersestion of $\mathcal{S}$ with the nilpotent cone is a surface
with an isolated double point singularity of type $A$, $D$, or $E$.
The exact correspondence is given by
\begin{align}
\ger{sl}_k &\longrightarrow A_{k-1} \\
\ger{sp}_{2k} &\longrightarrow D_{k+1} \\
\ger{so}_{2k} &\longrightarrow D_{k} \\
\ger{so}_{2k+1} &\longrightarrow A_{2k-1}
\end{align}
where the defining equations of these singularities are given by
\begin{align}
A_k &: \phantom{adad} Z^{k+1} + X^2 + Y^2 = 0\\
D_k &: \phantom{adad} Z^{k-1} + ZX^2 + Y^2 = 0. 
\end{align}
Moreover, the restriction of the adjoint quotient to $\mathcal{S}$
induces a semiuniversal deformation of the corresponding singularity.
\epp

\section{The Symplectic Case}\lb{sp}

\subsection{Nilpotent Orbits}\lb{sporbits}
 
Define linear functionals $e_i \in {\ger h}^*$ for $i = 1, \dots, m$ by
$e_i (x) = x_i$ where $x = \mbox{diag}(x_1, \dots, x_m) \in D_m \cong \ger h$.
The root system of $\ger{sp}_{2m}(\Cc)$ is
$\{ \pm e_i \pm e_j , \pm 2 e_i \mid 1 \leq i,j \leq n , i \neq j \}$.
We choose our positive roots to be  $\{ e_i \pm e_j , 2 e_k 
         \mid 1 \leq i < j \leq m , 1 \leq k \leq m \}$. 
Let $E_{i,j}$ denote the $2m \times 2m$ matrix having 1 in the
$(i,j)$th place and zeros everywhere else.
The $\alpha$-root space corresponding to a simple root $\alpha$ is the
root vector $X_{\alpha}$ given by
\begin{align}
X_{e_i-e_j} &= E_{i,j} - E_{j+m , i+m} \\
X_{e_i + e_j} &= E_{i, j+m} + E_{j, i+m} \\
X_{-e_i - e_j} &= E_{i+m , j} + E_{j+m , i} \\
X_{2e_i} &= E_{i,i+m} \\
X_{-2 e_i} &= E_{i+m , i}.
\end{align}  
For each 
$n = 0, 1, \dots, m$ let $X_n$  
be the matrix given by
\[
X_n = \sum_{i \neq m-n} X_{e_i - e_{i+1}}   
        + X_{2e_{m-n}}
        + X_{2e_{m}}.
\]
That is 
\[
X_n = 
\left[
\begin{tabular}{c c | c c}
  $J_{m-n}$ &        & $\mathcal{E}_{m-n,m-n}$    &           \\
            &  $J_n$ &              &  $\mathcal{E}_{n,n}$    \\ \hline
            &        & $-J_{m-n}^t$ &           \\
            &        &              & $-J_n^t$    
\end{tabular}
\right]
\]
where $J_l$ is the principle nilpotent Jordan block of size $l$ and
the matrix $\mathcal{E}_{j,k}$ is the $j \times k$ matrix that has 1 
in position $(j,k)$ and 
zeros everywhere else. Then $X_n$ is a representative
of the nilpotent orbit in $\ger{g}$ corresponding to the partition 
$[2(m-n), 2n]$ (\cite{CM} 5.1). We will assume that $n \leq m/2$,
or equivalently, $2(m-n) \geq 2n$, since if we replace $n$ with $m-n$ we do 
not change the orbit.
When $n=0$
we get a representative of the principle orbit. 
When $n=1$ we get a representative of the subregular orbit.
Also, $\mbox{codim}\mathcal{O}_{X_n} = m + 2n$ (\cite{CM} 6.1).

\subsection{Transverse Slices}\lb{spslice}

Now we describe a slice to the orbit of $X_n$.
Let $\{a_i, y_i, z_i, d_j \}$ be arbitrary complex numbers
for $1 \leq i \leq n$ and 
$1 \leq j \leq m-n$. Define a vector space $V_n$ by
\begin{align}
V_n = \sum_{i=1}^{m-n} d_{m-n-i+1} X_{-2e_i}
    &+ \sum_{i=m-n+1}^{m} a_{m-i+1} X_{-2e_i} \\
    &+ \sum_{j=1}^n y_j X_{-e_1 - e_{m-n+j}}
    + \sum_{j=1}^n z_j X_{-e_{j+1} - e_m}.
\end{align}
In order to make this clearer, we define an $m \times m$ symmetric matrix
$$
M(a_i, y_i, z_i, d_j)=
\left[
  \begin{tabular}{c c c c c | c c c c}
    $d_{m-n}$ & & & &      &   $y_1$ & $y_2$ & $\cdots$ & $y_n$ \\
      & $d_{m-n-1}$ & & &    &   & & & $z_1$ \\
      & & $\ddots$ & & &   & & & $\vdots$ \\
      & & & $\ddots$ & &   & & & $z_n$\\
      & & & & $d_{1}$&   & & &  \\ \hline
    $y_1$ & & & &   &                   $a_n$ & & & \\
    $y_2$ & & & &   &                     & $a_{n-1}$ & & \\
    $\vdots$ &&&&   &                  &  &  $\ddots$ & \\
    $y_n$ & $z_1$ & $\cdots$ & $z_n$ &  &  & && $a_1$
  \end{tabular}
\right].
$$
So for all $\{a_i, y_i, z_i, d_j \}$ we have
\[
   \left[
    \begin{tabular}{c | c}
       0 & 0 \\ \hline
    $M(a_i, y_i, z_i, d_j)$ & 0
    \end{tabular}
   \right]           \in V_n.
\]
We set
\[
\mathcal{S}_n = X_n + V_n
\]
and write elements of this set as
$S = S(a_i, y_i, z_i, d_j)$.
Then $\mathcal{S}_n$ is an affine subspace of $\ger g$. Notice that 
$V_n$ has a basis consisting
of root vectors $X_{\alpha}$ that correspond to negative simple roots.

A straightforward calculation 
shows that $\mathcal{S}_n$ is transverse to the orbit of
$X_n$. We will carry out this calculation in a later section.
For now we just state the result:

\bpp\lb{spmainprop}
Let $n \leq m/2$. The affine subspace $\mathcal{S}_n$ is 
a $\lambda$-invariant transverse slice to the 
adjoint orbit of $X_n$ at the point $X_n$. Moreover, 
for $S=S(a_i, y_i, z_i, d_j) \in \mathcal{S}_n$ the characteristic 
polynomial of $S$ is given by
\beq\lb{main}
\chi_S (t) = A(t)D(t) + (-1)^m B(t) B(-t)
\end{equation}
where
\begin{align}
A(t) &= t^{2n} - a_1 t^{2(n-1)} + \cdots
           + (-1)^{n-1}a_{n-1} t^2 + (-1)^{n} a_n \lb{spA}\\
D(t) &= t^{2(m-n)} - d_{1} t^{2(m-n-1)} + \cdots
           + (-1)^{m-n-1} d_{m-n-1} t^2 + (-1)^{m-n} d_{m-n} \lb{spD}\\
B(t) &= y_1 +y_2 t + \cdots + y_n t^{n-1}
         - t^n (z_1 - z_2 t + \cdots + (-1)^{n-1} z_n t^{n-1}). \lb{spB}
\end{align}
\epp

\brm 
Actually, our definition of $\mathcal{S}_n$ above only makes sense when $n$
is strictly less than $m/2$. When $n=m/2$ we have to modify our
definition slightly. This does not affect Proposition \ref{spmainprop}, however.
\erm

As is the previous section, for any $\tau \in \mbox{Sym}_{m}(\Cc)$
we set $\mathcal{Y}_{n,\tau} = \chi^{-1}(\tau) \cap \mathcal{S}_n$.
In terms of the coordinates $\{ a_i, y_i, z_i, d_j \}$ 
of $\mathcal{S}_n$, the variety $\mathcal{Y}_{n,\tau}$  
is described by a set of
$m$ algebraic equations that can be grouped into one:

\beq\lb{main}
A(t)D(t) + (-1)^m B(t) B(-t) = P_{\mu}(t)
\end{equation} 
where $P_{\mu}(t)$ is the unique monic polynomial with
roots given by $\mu$, and $\mu$ is the image of $\tau$ under
the embedding $\mbox{Sym}_m (\Cc) \into \mbox{Sym}_{2m} (\Cc)$. The
$m$ equations are given by equating
the coefficients of $t$ in \eqref{main}.

We define the following polynomials,
\beq\lb{UV}
U(t) = \frac{i^{m-1}}{2}(B(t)+B(-t)) \qquad \mbox{and} \qquad 
V(t) = \frac{i^{m}}{2t}(B(t)-B(-t)) 
\eeq
so that $U(t)^2 + t^2 V(t)^2 = (-1)^{m-1} B(t)B(-t)$.
We notice that both $U$ and $V$ are polynomials containing
only even powers of $t$. For any polynomial, say
$q(t)$, having only even powers of $t$, we can
define a new polynomial $\hat q (t)$ of half
the original degree by specifying that $\hat q (t^2) = q(t)$.
Since ${\hat P}_{\mu} (t) = P_{\tau} (t)$, we can rewrite \eqref{main}
to obtain
\beq\lb{main2}
P_{\tau}(t) + \hat U (t)^2 + t \hat V (t)^2
                 = \hat A (t) \hat D (t)  . 
\eeq

\begin{example}
Setting $n=1$ we see that in the sub-regular case, the varieties
we obtain have their defining equation given by
$P_{\tau}(t) + (-1)^{m-1} {y_1}^2 + t (-1)^m {z_1}^2
                 = (t- a_1) \hat D (t)$. Once we know that
$a_1$ is a root of $P_{\tau}(t) + (-1)^{m-1} {y_1}^2 + t (-1)^m {z_1}^2$,
we can recover $\hat D (t)$ uniquely. 
Thus, after a simple change of cooridinates we have
\beq\lb{spsubreg}
\mathcal{Y}_{1,\tau} = \{u,v,z \mid
    P_{\tau}(z) + u^2 + z v^2 =0 \}.
\eeq
Now, if we take $\tau = 0$ then the defining equation we obtain for
$\mathcal{Y}_{1,0}$ is 
\[
z^m + u^2 + z v^2 = 0.
\]
This is the type $D_{m+1}$
Kleinian singularity. Of course, setting $n=0$ is the same as 
intersecting the slice $\mathcal{S}_n$
with the nilpotent cone so we recover Slodowy's result (Proposition
\ref{Slodowy}).

We write $\Sigma_{\tau}  = \mathcal{Y}_{1,\tau}$. 
Because $\mathcal{S}_1$ is $\lambda$-invariant (which we prove next), 
the remarks in Section
\ref{geom} apply here to give us that $\Sigma_{\tau}$ is
smooth for $\tau \in R_{\ger{sp}_{2m}(\Cc)} = \mbox{Conf}_{m}^{\,*}(\Cc)$.
\end{example}

\subsection{$\lambda$-Invariance}\lb{inv}

Recall from Section \ref{geom} that a transverse slice 
$\mathcal S$ of a nilpotent
orbit $\mathcal{O}_X$ is said to be $\lambda$-invariant if
it is integral to the vector field 
\beq\lb{vectorfield}
\xi_Y = Y - \frac{1}{2}[H,Y]
\eeq
where $ \{ H , N^+ , N^- \} $ is some JM triple with $X = N^+$.
We show here that the slices $\mathcal{S}_n$ constructed above
are $\lambda$-invariant.

For any positive integer $k$, define two $k \times k$  matrices 
$l_k$ and $m_k$ by
\begin{align}
l_k &=
\left[
  \begin{tabular}{c c c c}
    $2k-1$ &    &    & \\
          & $2k-3$ & &   \\
       &   &  $\ddots$ &  \\
     & & & 1        
  \end{tabular}
\right]         \\
m_k &= 
\left[
  \begin{tabular}{c c c c c}
          0  &           &             &        &        \\
      $2k-1$ &  0        &             &        &        \\
             & $2(2k-2)$ &  0          &        &         \\
             &           &    $\ddots$ & $\ddots$    &         \\
             &           &             & $(k-1)(k+1)$ &    0    
  \end{tabular}
\right]
\end{align} 
Then we set
\begin{align}
H_n &=
\left[
  \begin{tabular}{c c | c c}
    $l_{m-n}$ &    &    & \\
          & $l_n$ & &   \\ \hline
       &   &  $- l_{m-n}$ &  \\
     & & & $-l_n$
  \end{tabular}
\right]               \\
N_n^- &=
\left[
  \begin{tabular}{c c | c c}
      $m_{m-n}$  &           &             &      \\
                         & $m_n$   &             & \\ \hline
          $(m-n)^2 \mathcal{E}_{m-n,m-n}$   &    &  $-m_{m-n}$  &  \\
             &  $ n^2 \mathcal{E}_{n,n}$         &     & $-m_n$  
  \end{tabular}
\right]
\end{align}
It is easy to check that $\{ H_n, X_n, N_n^- \}$ is a JM triple in
$\ger{sp}_{2m} (\Cc)$. It defines a vector field $\xi_n$ as in
equation \eqref{vectorfield} and an easy calculation shows
that $\mathcal{S}_n$ is integral to $\xi_n$. In fact, as we noted
above, $\mathcal{S}_n$ looks like $X_n$ plus an independent combination
of negative simple root vectors. The $\Cc^*$-action acts on these
root vectors by multiplication, hence $\mathcal{S}_n$ is invariant
under the action.

In terms of the coordinates of $\mathcal{S}_n$ the $\Cc^*$-action
is given by
\beq
\lambda_r (S (a_i, y_i, z_i, d_j)) =
       S ( r^{2i} a_i, r^{m-i+1} y_i, r^{m+n-i+1} z_i, r^{2j} d_j).
\eeq
Also, since $\lambda_r$ multiplies each eigenvalue by $r$, then on the level
of $\ger h /W \cong \mbox{Sym}_m (\Cc)$ it will multiply each 
$\tau \in \mbox{Sym}_m (\Cc)$ by $r^2$. This means that
\beq
\lambda_r (\mathcal{Y}_{n, \tau}) = \mathcal{Y}_{n, r^2 \tau}.
\eeq

Hence, the results of Section \ref{geom} apply to the slices
$\mathcal{S}_n$. Namely, 
the restriction to $\mbox{Conf}_m^{\,*} (\Cc)$ of
$\chi | \mathcal{S}_n$
is a differentiable fiber bundle. The fibers 
$\mathcal{Y}_{n, \tau}$ of this restriction 
are smooth 2n-dimensional manifolds.

\subsection{An Open Holomorphic Embedding}

Modulo $\hat A (t)$ the defining equation of $\mathcal{Y}_{n, \tau}$
in \eqref{main2} is formally identical to that of $\Sigma_{\tau}$
in \eqref{spsubreg}. As Manolescu does for the type-$A$ case, we
exploit this formality to construct an algebraic morphism
from $\mathcal{Y}_{n, \tau}$ in to the Hilbert scheme of
$n$ points on the complex surface $\Sigma_{\tau}$. The proof of the 
following theorem is identical to 
Manolescu's proof in the type-$A$ case. We include it here for
the sake of completeness. 

\btm\lb{sphilb}
Let $n \leq m/2$ and 
$\tau \in R_{\ger{sp}_{2m}}(\Cc) \subseteq \mbox{Conf}_{m}^{\:*}(\Cc)$. 
There is an open holomorphic embedding 
\beq
j: \mathcal{Y}_{n,\tau} \to \mbox{Hilb}^n (\Sigma_{\tau}),
\eeq
where $\Sigma_{\tau}$
is the affine surface in $\Cc^3$ described by equation \eqref{spsubreg}.
\etm

\begin{proof} We know that 
$\mathcal{Y}_{n, \tau} = \mbox{Spec} (R)$, where $R$
is the quotient of the polynomial ring in the $m+2n$
coefficients $a_i, y_i, z_i, d_j$ by the ideal
generated by the algebraic relations in \eqref{main2}. We
think of $\hat A (t)$, $\hat D (t)$, $\hat U (t)$, and $\hat V (t)$
as elements of $R[t]$. Define
$$
\mathcal{R} = R[u,v,z]/(P_{\tau}(z) +  u^2 + z v^2 ).
$$
Notice that 
$\Sigma_{\tau} = 
   \mbox{Spec} (\Cc [u,v,z]/(P_{\tau}(z) +  u^2 + z v^2 ))$
so that 
$$
\mbox{Spec}( \mathcal{R} )= \mathcal{Y}_{n, \tau} \times \Sigma_{\tau}.
$$
Consider the map $\psi : \mathcal{R} \to R[t]/(\hat A (t))$
defined by 
$$
\psi (Q(u,v,z))= Q(\hat U (t), \hat V (t), t).
$$
It is clear that this map is 
well defined and surjective. Let $\mathcal{K}$ be
its kernel. Then $\mathcal{R}/\mathcal{K}$ is isomorphic
to $R[t]/(\hat A (t))$. But $ R[t]/(\hat A (t)) \cong R^n$ since
$\hat A (t)$ is a monic polynomial of degree $n$
in $R[t]$. We define the closed
subscheme 
$$
Z= \mbox{Spec}(\mathcal{R}/\mathcal{K}) \subseteq
        \mbox{Spec} ( \mathcal{R} )= \mathcal{Y}_{n, \tau} \times \Sigma_{\tau}.
$$
Since $\mathcal{R}/\mathcal{K}$ is a free $n$-dimensional module
over $R$, it follows that the projection
$$
Z \subseteq \mathcal{Y}_{n, \tau} \times \Sigma_{\tau}
          \to \mathcal{Y}_{n, \tau} = \mbox{Spec}( R)
$$
exhibits $Z$ as a flat family of 0-dimensional subschemes of $\Sigma_{\tau}$
of length $n$. This defines a map
$$
j : \mathcal{Y}_{n, \tau} \to \mbox{Hilb}^n (\Sigma_{\tau}).
$$

Now that we have defined the map we need to prove that it is 
an open embedding. 

With its reduced scheme structure, the points of $\mathcal{Y}_{n,\tau}$
are 4-tuples of polynomials $(\hat A, \hat D, \hat U, \hat V )$ 
in $\Cc [t]$ satisfiying the defining equation in \eqref{main2}. The points
in $\mbox{Hilb}^n (\Sigma_{\tau})$ can be identified with ideals
$\mathcal I$ in the coordinate ring
$R_{\Sigma} = \Cc [u,v,z]/(P_{\tau}(z) + u^2 + z v^2 )$ of $\Sigma_{\tau}$
such that
$\dim_{\Cc} (R_{\Sigma} / \mathcal I ) = n$. Explicitly,
the morphism $j$ is given by 
\beq
j ( \hat A, \hat D, \hat U, \hat V) =
   \{ Q(u,v,z) \mid \hat A (t) 
           \mbox{ divides } Q( \hat U (t), \hat V (t), t) \}.
\eeq
Let $R_0 = \Cc [z]$. This is a subring of $\Cc [u,v,z]$, so
$R_1 = R_0 / (R_0 \cap ( P_{\tau} (z) + u^2 + z v^2 )) \cong \Cc [z]$ 
is a subring of $R_{\Sigma}$.

Now, $\mbox{Hilb}^n (\Sigma_{\tau})$  is irreducible and has the same 
dimension as $ \mathcal{Y}_{n, \tau} $, so in order to prove that $j$
is an open embedding it suffices to show that it is injective.

So pick $( \hat A_i, \hat D_i, \hat U_i, \hat V_i)$ for $i=1,2$ that
map to the same ideal $\mathcal I$ under $j$. Then
\[ 
\mathcal{I} \cap R_1 = \{ Q \in \Cc [z] \mid \hat A_i (z) 
        \mbox{ divides } Q(z) \},
\]
where $i$ can be either 1 or 2. Thus, we must have that
$\hat A_1 (z)$ divides $\hat A_2 (z)$. But each $\hat A_i$ is 
a monic polynomial of degree $n$, so we must have $\hat A_1 = \hat A_2$.

Next, note that $u - \hat U_1 (z)$ and $u - \hat U_2 (z)$ are
in $\mathcal{I}$ so 
$\hat U_1 (z) - \hat U_2 (z) \in \mathcal{I} \cap R_1$. But
$\hat U_1 (z) - \hat U_2 (z)$ has degree
at most $n-1 < n = \deg \hat A$, so $\hat A$ dividing 
$\hat U_1 - \hat U_2$ implies $\hat U_1 = \hat U_2$. Similarly,
we must have $\hat V_1 = \hat V_2$. Also, the relation \eqref{main2}
determines $\hat D$ uniquely in terms of $\hat A$, $\hat U$, $\hat V$, and 
$P_{\tau}$. Thus, we must have $\hat D_1 = \hat D_2$ so that
the map $j$ is injective.  
\end{proof}

\brm
Note that everything in the above proof holds for arbitrary
$\tau \in \ger h / W$. That is, we can define an algebraic map
map $j : \mathcal{Y}_{n,\tau} \to \mbox{Hilb}^n (S_{\tau})$
and show it is injective. However, $S_{\tau}$ may have singularites
for general $\tau$.
\erm

Let $\mathcal I \subseteq R_{\Sigma}$ be an ideal describing a subscheme
$Z = \mbox{Spec} R_{\Sigma} / \mathcal I $ in $\mbox{Hilb}^n (\Sigma_{\tau})$.
Then $\mathcal I \cap R_1 $ corresponds to a subscheme of $\Cc$, namely the
image of $Z$ under the map $i : \Sigma_{\tau} \to \Cc$ defined by
\[
i (u,v,z) = z.
\]
Since $R_1 / (\mathcal I \cap R_1 )$ injects into $R_{\Sigma} / \mathcal I$, we
must have that $i(Z)$ has length at most $n$.

\bpp\lb{spimage}
The image of $j$ in $\mbox{Hilb}^n (\Sigma_{\tau})$ consists of all subschemes
$Z$ such that $i (Z)$ is a subscheme of $\Cc$ of length exactly $n$. Moreover,
given a point $(\hat A,\hat D,\hat U,\hat V)$ in $\mathcal{Y}_{n,\tau}$
its image under $\pi_{HC} \circ j$ is the unordered set of
$n$ points $(u_k, v_k, z_k) \in \mbox{Sym}_n (S_{\tau})$ where
$z_k$ are the roots of $\hat A (t)$, $u_k = \hat U (z_k)$, and 
$v_k = \hat V (z_k)$.
\epp 

The proof of Proposition \ref{spimage} is verbatim the same as in the 
type-$A$ case so we refer readers to \cite{Man} for the details.

\subsection{Proof of Proposition \ref{spmainprop}}

We start by assuming $n < m/2$ 
(the case $n \leq m/2$ will be handled separately). We let
$X_n$  be the representative of the nilpotent orbit
corresponding to the partition $[2(m-n),2n]$ that was given in 
Section \ref{spslice}.

Recall that the matrix $M_n (a_i, y_i, z_i, d_j)$ 
is defined by
$$
M_n (a_i, y_i, z_i, d_j)=
\left[
  \begin{tabular}{c c c c c | c c c c}
    $d_{m-n}$ & & & &      &   $y_1$ & $y_2$ & $\cdots$ & $y_n$ \\
      & $d_{m-n-1}$ & & &    &   & & & $z_1$ \\
      & & $\ddots$ & & &   & & & $\vdots$ \\
      & & & $\ddots$ & &   & & & $z_n$\\
      & & & & $d_{1}$&   & & &  \\ \hline
    $y_1$ & & & &   &                   $a_n$ & & & \\
    $y_2$ & & & &   &                     & $a_{n-1}$ & & \\
    $\vdots$ &&&&   &                  &  &  $\ddots$ & \\
    $y_n$ & $z_1$ & $\cdots$ & $z_n$ &  &  & && $a_1$
  \end{tabular}
\right]
$$
and that for any $\{a_i, y_i, z_i, d_j \}$ we defined the
element $S=S(a_i, y_i, z_i, d_j) \in \ger g$ by
$$
S(a_i, y_i, z_i, d_j) =  X_n +
   \left[
    \begin{tabular}{c | c}
       0 & 0 \\ \hline
    $M_n (a_i, y_i, z_i, d_j)$ & 0
    \end{tabular}
   \right].
$$

We wish to show that $\mathcal{S}_n$ is transverse to 
$\mathcal{O}_{X_n}$.
The second part of Proposition \ref{spmainprop}, namely 
that $\chi_{S} (t) = A(t)D(t) + (-1)^m B(t) B(-t)$ for 
$A(t)$, $D(t)$,
and $B(t)$ given by equations \eqref{spA}-\eqref{spB}, is better
left to the reader. It is a straightforward computation.

Since $\dim \mathcal{S}_n = m +2n = \mbox{codim }[X_n, \ger{g}]$ 
(\cite{CM} 6.1), then  to show that $\mathcal{S}_n$ 
is transverse to the orbit through $X_n$ amounts to
showing that $V_n$ intersects $[X_n, \ger{g}]$ trivially.

Let $T=[t_{i,j}]$ and $C = [c_{i,j}]$ 
be arbitrary $n \times m$ matrices.
We define operations
$a$ and $b$ as follows:
\begin{align}
a(C) &= \left[
  \begin{tabular}{c | c c c}
   0            & $c_{1,1}$  & $\cdots$ & $c_{1,m-1}$\\ \hline
   $c_{1,1}$    &            &          &            \\
   $\vdots$     &            & $c_{i-1,j} + c_{i,j-1}$ & \\
   $c_{n-1,1}$  &            &                         &
  \end{tabular} 
 \right] \\
b(T,C) &= \left[
  \begin{tabular}{c | c c c}
  $t_{2,1}$ &         &                        & \\
  $\vdots$  &         & $t_{i+1,j}-t_{i,j-1}$ & \\
  $t_{n,1}$ &         &                       &\\ \hline
  $c_{n,1}$ &         &  $c_{n,j} - t_{n, j-1}$ &
\end{tabular}\right].
\end{align}

\blm\lb{lemmaA}
Let $T=[t_{i,j}]$ and $C = [c_{i,j}]$ be $n \times n$ matrices,
with $C$ symmetric, and suppose 
that  $b(T,C)=0$ and $a(C)=\mbox{diag}(x_1, \dots, x_n)$. Then
$x_1 = \cdots = x_n = 0$.
\elm
\begin{proof} Looking at
the $(n-k, j-k)$ entry of $b(T,C)$ for $k= 0,1,\dots, j-1$ we see that
if $b(T,C)=0$ then
$$
c_{n,j} = t_{n,j-1}=t_{n-1,j-2}=\cdots =t_{n-j+2,1}=0
$$
for all $j$. Since $C$ is also symmetric we have
$$
a(C) = \left[
  \begin{tabular}{c | c c c | c}
   0            & $c_{1,1}$  & $\cdots$ & $c_{1,n-2}$ &  $c_{1,n-1}$\\ \hline
   $c_{1,1}$    &            &          &             &  $c_{2,n-1}$ \\
   $\vdots$     &            & $c_{i-1,j} + c_{i,j-1}$& & $\vdots$ \\
   $c_{n-2,1}$  &            &                   & & $c_{n-1, n-1}$ \\ \hline
   $c_{n-1,1}$  & $c_{n-1,2}$& $\cdots$ & $c_{n-1,n-1}$ & 0
  \end{tabular}
 \right].
$$
So if $a(C)=\mbox{diag}(x_1, \dots, x_n)$ then we have
$$
x_i = c_{i-1,i} + c_{i,i-1} = 2c_{i,i-1}.
$$
But the $(i+k,i-k)$th entry for $k=1,2,\dots, l=\min \{i-1, n-i \}$ gives us
\[
c_{i,i-1}= -c_{i+1,i-2} = c_{i+2,i-3} = \cdots = (-1)^{l-1}c_{i+l-1,i-l}= 0.\qedhere
\]
\end{proof}

Let us denote by $K_{y_i,z_j}$
the $n \times m$ matrix $\{k_{i,j} \}$ defined by 

\begin{align}
&&k_{i,1} &= y_i &&\mbox{for } i=1,2,\dots,n &&\\ 
&&k_{n,j} &= z_{j-1} &&\mbox{for } j=2,3,\dots, n+1 &&\\
&&k_{i,j} &= 0 &&\mbox{otherwise}.&&
\end{align}

\blm\lb{lemmaB}
Let $T=[t_{i,j}]$ and $C=[c_{i,j}]$ be arbitrary $n \times m$ matrices,
let $R=[r_{i,j}]$ be an arbitrary $m \times n$ matrix and suppose that
$b(T,C)=0$, $b(R,C^{\top})=0$, and $a(C)=K$. Then
we have $K=0$.
\elm
\begin{proof} Let $j \leq n$ and $k=0,1,\dots, j-1$. From the
$(n-k, j-k)$th entry of $b(T,C)$ we see, as in the previous theorem, that
if $b(T,C)=0$ then
$$
c_{n,j} = t_{n,j-1}=t_{n-1,j-2}=\cdots =t_{n-j+2,1}=0.
$$
That is, $c_{n,1}= c_{n,2}=\cdots =c_{n,n}=0$.
Similarly, from $b(R,C^{\top})=0$ we see that
$$
c_{1,m} = c_{2,m} = \cdots = c_{n,m}=0.
$$
Thus, we can write
$$
a(C)=  \left[
  \begin{tabular}{c | c c c c | c}
   0  & $c_{1,1}$  & $\cdots$ & $c_{1,n}$ & $\cdots$ &  $c_{1,m-1}$\\ \hline
   $c_{1,1}$    & &          &          &             &  $c_{2,m-1}$ \\
   $\vdots$     &     & $c_{i-1,j} + c_{i,j-1}$& & & $\vdots$ \\
   $c_{n-2,1}$  &            &                 & & & $c_{n-1, m-1}$ \\ \hline
   $c_{n-1,1}$  & $c_{n-1,2}$& $\cdots$ & $c_{n-1,n+1}$ & $ \ast \ast \ast$ & $c_{n,m-1}$
  \end{tabular}
 \right].
$$
So if we have $a(C) = K$ then we have
$$
y_i = c_{i-1,1}= -c_{i-2, 2} = \cdots = (-1)^{i}c_{1,i-1}=0.
$$
Also, let $j \in \{1, 2, \dots, n\}$ and set $l = \mbox{min} \{n-1, m-1-j\}$.
Then
\[
z_j = c_{n-1,j+1} = -c_{n-2,j+2} = \cdots = (-1)c_{n-l,j+l} =0. \qedhere
\]
\end{proof}

So now let us write a $2m \times 2m$ matrix $T$ as four blocks of size $m$, 
each of which is then futher broken down into blocks:
\beq
T = \left[
\begin{tabular}{c | c}
  \begin{tabular}{c c}
    $T_{1,1}$ & $T_{1,2}$ \\
    $T_{2,1}$ & $T_{2,2}$
  \end{tabular}    
      &
\begin{tabular}{c c}
    $B_{1,1}$ & $B_{2,1}^{\top}$ \\
    $B_{2,1}$ & $B_{2,2}$
  \end{tabular}
       \\\hline
  \begin{tabular}{c c}
    $C_{1,1}$ & $C_{2,1}^{\top}$ \\
    $C_{2,1}$ & $C_{2,2}$
  \end{tabular}
        &
 \begin{tabular}{c c}
    $-T_{1,1}^{\top}$ & $-T_{2,1}^{\top}$ \\
    $-T_{1,2}^{\top}$ & $-T_{2,2}^{\top}$
  \end{tabular}
\end{tabular} \right].
\eeq

We want, for instance, $T_{1,1}$ and $T_{2,2}$ to be of size $m-n$
and $n$, respectively. So then $T \in \ger{sp}_{2m}(\Cc)$ 
if $B_{k,k}$ and $C_{k,k}$ are symmetric for $k=1,2$. So assume 
$T \in \ger{sp}_{2m}(\Cc)$. Then 
$$
[X_n ,T] =  \left[
\begin{tabular}{c | c}
  \begin{tabular}{c c}
    $b(T_{1,1},C_{1,1})$ & $b(T_{1,2},C_{2,1}^{\top})$ \\
    $b(T_{2,1},C_{2,1})$ & $b(T_{2,2},C_{2,2})$
  \end{tabular}
      &
        \Large{$\ast$}
       \\\hline
  \begin{tabular}{c c}
    $a(C_{1,1})$ & $a(C_{2,1}^{\top})$ \\
    $a(C_{2,1})$ & $a(C_{2,2})$
  \end{tabular}
        & \Large{$\ast$}
\end{tabular} \right].
$$
So if $S(a_i, y_i, z_i, d_j)-X_n =[X_n ,T]$ then we have 
\begin{align}
b(T_{1,1}, C_{1,1}) &= 0, & a(C_{1,1}) 
           &= \mbox{diag}(d_1, \dots, d_{m-n}),\lb{trans1}\\
b(T_{2,2},C_{2,2}) &= 0, &  a(C_{2,2}) 
           &= \mbox{diag}(a_1, \dots, a_{n}),\lb{trans2}
\end{align}
and
\begin{align}
b(T_{2,1},C_{2,1}) &= 0, 
          & b(T_{1,2},C_{2,1}^{\top}) &=0, 
          & a(C_{2,1}) &= K{y_i, z_j}.\lb{trans3}
\end{align}

Thus, applying Lemma \ref{lemmaA} to \eqref{trans1} and \eqref{trans2}, and 
Lemma \ref{lemmaB} to \eqref{trans3}, we see that 
$\mathcal{S}_n - X_n$ 
does indeed intersect $[ X_n , \ger g ]$ trivially.

Now we turn to the case when $n=m/2$. Of course we need
$m$ to be even. In this case we have to modify our 
definition of $S(a_i, y_i, z_i, d_j)$
slightly. Namely, we set
\[
S(a_i, y_i, z_i, d_j) =  X_n +
   \left[
    \begin{tabular}{c | c}
  0   & 0  \\ \hline
    $M_n^{'} (a_i, y_i, z_i, d_j)$ & 0
    \end{tabular}
   \right]   + z_n X_{e_m - e_n}
\]
where
\[
M_n^{'} (a_i, y_i, z_i, d_i)=
\left[
  \begin{tabular}{c c c  c | c c c c}
    $d_{n}$ & & &       &   $y_1$ & $y_2$ & $\cdots$ & $y_n$ \\
      & $d_{n-1}$ & &     &   & & & $z_1$ \\
      & & $\ddots$ & &    & & & $\vdots$ \\
      & & &      $d_1 + z_n^2$ &    & & & $z_{n-1}$\\ \hline
    $y_1$ & & & &                      $a_n$ & & & \\
    $y_2$ & & & &                        & $a_{n-1}$ & & \\
    $\vdots$ &&&&                     &  &  $\ddots$ & \\
    $y_n$ & $z_1$ & $\cdots$ & $z_{n-1}$   &  & && $a_1$
  \end{tabular}
\right].
\]
So, in light of what we have already proven, to show that this slice is transverse
to the orbit of $X_n$ all we need to show is that  if
$S(a_i, y_i, z_i, d_j)-X_n =[X_n ,T]$ then $z_n$
must equal 0. But this is clear from
Lemma \ref{lemmaB} since the equation $b(T_{1,2} , C^{\top}_{2,1} ) = 0$
gives us $c_{n,n} = 0$ and the equation 
$b( T_{2,1}, C_{2,1}) = z_n \mathcal{E}_{n,n}$ gives us $z_n = c_{n,n}$.

\section{The Orthogonal Case}

We start the orthogonal case with the type-$D$ algebras, $\ger{so}_{2m}(\Cc)$.
Once we have all the results in for this case we can hopefully make quick work 
of the type-$B$ case by viewing $\ger{so}_{2m+1}(\Cc)$ as a 
trivial $2m$-bundle over $\ger{so}_{2m}(\Cc)$. 
Constucting the transverse slices is somewhat more 
complicated for the orthogonal algebras 
than for $\ger{sl}_{2m}$ and $\ger{sp}_{2m}$. 

\subsection{Nilpotent Orbits in $\ger{so}_{2m} (\Cc)$}

Let $\ger g = \ger{so}_{2m}(\Cc)$. We take the standard Cartan 
subalgebra $\ger h \subset \ger g$, the subspace in $\ger g$ of all 
diagonal matrices. This is
identical to the Cartan subalgebra in the symplectic case.
The root system in this case is
$\{ \pm e_i \pm e_j \mid 1 \leq i,j \leq m , i\neq j \}$ and we
take $\{ e_i \pm e_j \mid 1 \leq i < j \leq m \}$ to be the set of
positive roots. The $\alpha$-root space is spanned by the vector $X_{\alpha}$
given by

\begin{align}
X_{e_i - e_j} &= E_{i,j} - E_{m+j,m+i} \\
X_{e_i + e_j} &= E_{i, m+j} - E_{j,m+i} \phantom{ajaj} (i < j) \\
X_{- e_i - e_j} &= E_{m+i, j} - E_{m+j,i} \phantom{ajaj} (i < j).
\end{align}

For each $n=0,1, \dots , m$ let $X_n$ be the matrix given by
\beq
X_n = \sum_{i \neq m-n} X_{e_i - e_{i+i}} 
              + X_{e_{m-n-1} + e_{m-n}}
              + X_{e_{m-n} + e_{m}}.
\eeq
That is,
\beq\lb{Xsoeven}
X_n =
\left[
\begin{tabular}{c c | c c}
   $J_{m-n}$ &       &  $\mathcal{F}_{m-n}$          & $\mathcal{E}_{m-n,n}$ \\ 
             & $J_n$ &  $-{\mathcal E}_{n,m-n}$ &    \\ \hline
             &       &   $-J_{m-n}^{\top}$     &       \\ 
             &       &                         & $-J_n^{\top}$
  \end{tabular}
\right]
\eeq
where $J_l$ is the principle $l \times l$ nilpotent Jordan block, 
$\mathcal{F}_l$ is the $l \times l$ matrix having 
1 in the $(l-1,l)$th place and -1 in the $(l, l-1)$th place,
and as before ${\mathcal E}_{j,k}$ is the 
$j \times k$ matrix having 1 in the $(j,k)$th place.

The matrix $X_n$ is a representative of the  nilpotent orbit in 
$\ger g$ corresponding to the  partition $[2(m-n)-1  ,2n+1]$. In order to 
avoid repitions we stipulate
that $2n+1 \leq m$. For these $n$, the codimension of $\mathcal{O}_{X_n}$ is
$m+2n$. Setting $n=0$ gives us a representative of the principle orbit. Setting
$n=1$ gives us a representative of the subregular orbit.

\subsection{Transverse Slices}

Let $\{ a_i , y_i , z_j , d_k \mid 1 \leq i \leq n , 1 \leq j \leq n+1,
                                         1 \leq k \leq m-n-1 \}$
be a set of $m+2n$ complex coordinates. We define a transverse slice
$\mathcal{S}_n$ in terms of root vectors as follows:
\beq\lb{soslice1}
\begin{aligned}
S ( a_i , y_i , z_j , d_k ) =
           X_n &+ a_1 X_{e_{m-n}-e_m} - d_{1} X_{e_{m-n}-e_{m-n-1}} \\
 &+  \sum_{i=1}^{n-1} a_{n-i+1} X_{-e_{m-n+i}-e_{m-n+i+1}}
               +  \sum_{i=1}^{n} y_i X_{-e_{1}-e_{m-n+i}}   \\
 &-  \sum_{j=1}^{n+1} z_j X_{e_{m-n}-e_{j}}
               +  \sum_{k=1}^{m-n-1} d_{m-n-k} X_{-e_{k}-e_{k+1}}.
\end{aligned}
\eeq

This is not really as bad as it looks. If we write out equation
\eqref{soslice1} as a matrix we can see what is going on. First let us 
write two square antisymmetric matrices as follows:
\beq
M_d = \left[
\begin{tabular}{c c c c}
0  &  $d_{m-n-1}$  &  &  \\
$-d_{m-n-1}$ & 0 & $\ddots$ &   \\
   & $\ddots$ & $\ddots$ & $d_{1}$ \\
  &  &  $-d_{1}$ & 0
\end{tabular}
\right]
\eeq
\beq
M_a = \left[
\begin{tabular}{c c c c}
0  &  $a_n$  &  &  \\
$-a_n$ & 0 & $\ddots$ &   \\
   & $\ddots$ & $\ddots$ & $a_{2}$ \\
  &  &  $-a_{2}$ & 0
\end{tabular}
\right].
\eeq
Then write an $(m-n) \times n$ matrix $M_y$ and
an $(m-n) \times (m-n)$ matrix $M_z$ as follows:
\begin{align}
M_y &= \left[
\begin{tabular}{c c c}
$y_1$  &  $\dots$  & $y_n$ \\
    & 0 &    \\
    &  & 
\end{tabular}
\right]
&
M_z &= \left[
\begin{tabular}{c c c c c}
0 & &  &  &  \\
 & &&&\\
&& $\ddots$  &&\\
       &&&&\\
 $z_1$ & $\dots$ &  $z_{n+1}$ & & 0
\end{tabular}
\right].
\end{align}
Notice that $M_z$ is well defined since $2n+1 \leq m$,
or rather $n+1 \leq m-n$.

Lastly, if we let $E_{m-n,m-n-1}$ be the $(m-n) \times (m-n)$ matrix
with 1 in the $(m-n,m-n-1)$th place and zeros everywhere else, and then 
write
\[
M_{z,d} = M_z + d_{1} E_{m-n,m-n-1}
\]
then we can finally write \eqref{soslice1} as
\beq\lb{soslice}
S ( a_i , y_i , z_j , d_k ) = X_n +
\left[
\begin{tabular}{c c | c c}
$- M_{z,d}$ & $ a_1 \mathcal{E}_{m-n,n}$ & & \\
0  & 0 & & \\ \hline
$M_d$ & $M_y$ & $M_{z,d}^{\top}$ & 0 \\
$-M_y^{\top}$ & $M_a$ &   $-a_1 \mathcal{E}_{n,m-n}$ & 0
\end{tabular}
\right]
\eeq

Let $\mathcal{S}_n$ be the collection of all elements 
$S ( a_i , y_i , z_j , d_k ) \in \ger g$. Then $\mathcal{S}_n$ is an affine
subspace containing $X_n$.

\btm
Let $2n+1 \leq m$. 
The affine subspace $\mathcal{S}_n$ is $\lambda$-invariant and 
transverse to the adjoint orbit of $X_n$. Moreover, the characteristic
polynomial of an element $S = S ( a_i , y_i , z_j , d_k )$ in 
$\mathcal{S}_n$ is given by
\beq\lb{somain}
\chi_S (t) = t^2 A(t) D(t) + (-1)^m B(t) B(-t) 
\eeq 
where
\begin{align}
A(t) &= t^{2n} - a_1 t^{2(n-1)} + \cdots + (-1)^n a_n \\ \lb{soA}
D(t) &= 4 d_1 t^{2(m-n-2)} - 4 d_2 t^{2(m-n-3)} 
          + \cdots + (-1)^{m-n} 4 d_{m-n-1} \\ \lb{soD}
B(t) &= y_1 + y_2 t + \cdots + y_n t^{n-1}   \\ \lb{soB}
     & \phantom{12345} - t^n ( z_1 - z_2 t + 
        \cdots + (-1)^n z_{n+1} t^n + (-1)^{m-n} t^{m-n}). \notag
\end{align}
\etm
\begin{proof}
We will prove transversality in a later section. The characteristic polynomial
can be computed easily using elementary row operations as in the 
symplectic case. For $\lambda$-invariance, notice that since $X_n$ is
a sum of positive root vectors  
we can find a JM-triple $H_n, X_n, N_n^-$
where $H_n$ is diagonal. To do this explicitly,
define $k \times k$ matrices
$\alpha_k$ and $\beta_k$ as follows
\begin{align}
\alpha_k &= \left[
  \begin{tabular}{c c c c}
    $2(k-1)$ &&&\\
    & $2(k-2)$ &&\\
    && $\ddots$&\\
    &&& 0
  \end{tabular}
\right]   \\
\beta_k &= \left[
  \begin{tabular}{c c c c}
    $2k$ &&&\\
    & $2(k-1)$ &&\\
    && $\ddots$&\\
    &&& 2
  \end{tabular}
\right]  .
\end{align}
Then set 
\beq
H_n = \left[
  \begin{tabular}{c c | c c}
    $\alpha_{m-n}$ &&&\\
    & $\beta_n$ &&\\ \hline
    && $-\alpha_{m-n}$&\\
    &&& $-\beta_n$
  \end{tabular}
\right].
\eeq
The matrix $H_n$ will then be the semisimple part of a JM-triple
having $X_n$ as its nilpositive part. Clearly, $ad(H_n)$ preserves
the $\alpha$-root space for any principle root $\alpha$. The only question 
that arises is does $ad(H_n)$ preserve
the one dimensional space  spanned by
$X_{-e_{m-n-1} - e_{m-n}} - X_{e_{m-n} - e_{m-n-1}}$, because this
is the space that carries the $d_1$-cooridinate in $\mathcal{S}_n$.
But calculating the $ad(H_n)$ action is easy and
we see that it acts on this one dimensional space by multiplication
by -2.
\end{proof}

Recall from Section \ref{geom} that we have a 
$\Z / 2$-action on $\ger h / W$ that multiplies a single coordinate
by -1 and that this groups action defines the fibers
of the map $\epsilon : \ger h / W \to \mbox{Sym}_m (\Cc)$. 
Let us denote the $\Z / 2$-action by $c$. Let $\tilde{\tau}$ be any element
of $\ger h$ so we have
\beq
\tilde{\tau} = \left[
\begin{tabular}{c c}
$D$ & 0 \\
0 & $-D$
\end{tabular}
\right]
\eeq
for some $m \times m$ diagonal matrix $D= \mbox{diag}(x_1, \dots, x_m)$.
Define $p : \ger h \to \Cc$ by $p(\tilde{\tau}) = \Pi_{i=1}^{m} x_i$. Notice
that $p$ factors through the projection $\ger h \to \ger h / W$.
So for any $\tau \in \ger h /W$ we have
\begin{align}
\epsilon( c \tau) &= \epsilon( \tau) & &\mbox{and} & p(c \tau)&=-p(\tau).
\end{align}
Define $U(t)$ and $V(t)$ as in \eqref{UV}. Then \eqref{somain}
can be rewritten to give 
\beq\lb{soeqn}
\hat \chi_{S} (t) + \hat U (t)^2 + t \hat V (t)^2 
        = t \hat A (t) \hat D (t).
\eeq
Notice that the only $\hat U$ and $\chi_S$ have constant terms in \eqref{soeqn}.
Let us say that $-y$ is the constant term of $\hat U$ 
($y$ is equal to some constant multiple of $y_1$) 
and let us define $\hat W$ by 
\beq
\hat U (t) = t \hat W (t) - y.
\eeq 
Then \eqref{soeqn} can once again be rewritten as
\beq\lb{soeqn2}
\hat \chi_{S} (t) - y^2 - 2yt \hat W (t) 
      + t^2 \hat{W} (t)^2+ t \hat V (t)^2 = t \hat A (t) \hat D (t).
\eeq

Now let $\tau$ be an element of $\ger h /W$ and, as usual,
let $\mathcal{Y}_{n,\tau} = \chi^{-1}(\tau) \cap \mathcal{S}_n$.
The set of all elements of $\ger g$ having reduced
characteristic polynomial
equal to  $P_{\epsilon(\tau)} (t)$ is exaclty
the union $\mathcal{Y}_{n,\tau} \cup \mathcal{Y}_{n, c \tau}$.
So we can describe this union
as the set of all 
$\{ a_i, y_i, z_j, d_k  \}$ that satisfy
the equation
\beq\lb{soeqn3}
P_{\epsilon(\tau)} (t) - y^2 - 2yt \hat W (t)
      + t^2 \hat{W} (t)^2+ t \hat V (t)^2 =  t \hat A (t) \hat D (t).
\eeq
By equating the constant terms in \eqref{soeqn3} we obtain
\[
p(\tau)^2 = y^2
\] 
to which  we have the two solutions
\begin{align}
y &= p(\tau) & y &= -p(\tau) = p(c \tau).
\end{align}

For any $\tau \in \ger h / W$ let us define a polynomial 
$Q_{\tau}(t)$ of degree $m-1$ by
\beq
Q_{\tau}(t) = \frac{1}{t}(P_{\epsilon(\tau)} - p(\tau)^2).
\eeq 

We are now in a position  to give a description 
of $\mathcal{Y}_{n,\tau}$ as a variety: it is
equal to the set of all 
\beq
\{ a_i, y_j, z_k, d_l  \mid 1 \leq i \leq n, 
                            2 \leq j \leq n,
                            1 \leq k \leq n+1,
                            1 \leq l \leq m-n-1 \}
\eeq
that satisfy
the equation
\beq\lb{soeqn4}
Q_{\tau}(t) + \hat V (t)^2 + t \hat W (t)^2 - 2 p(\tau)  \hat W (t)
             = \hat A (t) \hat D (t).
\eeq

If  $\tau$ belongs to  $(\ger h / W)^0$ then we can simplify
this expression considerably. Recall that
this $\tau \in (\ger h / W)^0$ occurs when  $\tau$
has a zero eigenvalue or, equivalently, when $c \tau = \tau$. In this case,
$\epsilon$ is a bijection so we can write
$\tau = \{\tau',0\} \in \mbox{Sym}_{m}^{\,0}(\Cc)$
for some uniquely defined $\tau' \in \mbox{Sym}_{m-1}(\Cc)$. Also, we have 
$p(\tau) = 0$ and $Q_{\tau}(t) = P_{\tau'}(t)$ in this case.

Thus, for $\tau \in (\ger h / W)^0$ we see that the defining equation 
for $\mathcal{Y}_{n,\tau}$ is
\beq\lb{soeqn5}
P_{\tau'}(t) + \hat V (t)^2 + t \hat W (t)^2 
             =   \hat A (t) \hat D (t).
\eeq

\begin{example}
Let us consider the subregular case. This is when $n=1$.
After changing a few signs we see that 
\begin{align}
\hat V (t) &= z_1 + \delta_{m-1} t^{(m-1)/2 -1}\\
\hat W (t) &= z_2 - \delta_m  t^{m/2 -1}\\
\hat A (t) &= t- a_1
\end{align}
where $\delta : \Z \to \{0,1\}$ is the function that is 1 on even integers
and 0 on odd integers. Now once we know that the $a_1$ is a root of
the left-hand side of \eqref{soeqn4} we can recover $\hat D$ uniquely
in terms of $a_1$, $z_1$, $z_2$, and $\tau$. We put all this together
with a change of cooridnates of the form
$w=z_i \pm a_1^{k}$ where
$i$ is either  1, or 2 depending on if $\delta_m = 1,0$, repectively. This
gives us the following:
\beq\lb{sosubreg}
\mathcal{Y}_{1,\tau}= 
 \{ v,w,z \mid Q_{\tau}(z) + v^2 + z w^2 = p(\tau) w \}.
\eeq
This is a complex surface in $\Cc^3$ and we know from Section \ref{geom}
that it is smooth for $\tau \in R_{\ger{so}_{2m}}(\Cc)$.
Let us give a name to the varitey in \eqref{sosubreg} by writing
$\Gamma_{\tau} = \{ v,w,z \mid Q_{\tau}(z) + v^2 + z w^2 = p(\tau) w \}$.

If $\tau = \{\tau',0\} \in (\ger h / W)^0 = \mbox{Sym}_{m-1}(\Cc)$, then
\eqref{sosubreg} reduces to
\beq\lb{sosubreg2}
\mathcal{Y}_{1,\tau}=
 \{ v,w,z \mid P_{\tau'}(z) + v^2 + z w^2 = 0 \}.
\eeq
Hence, in this case we obtain the same surface as in the type-$C$ case
so we will denote it the same way and write
$\Gamma_{\{\tau',0\}} = \Sigma_{\tau'}$.

If we take $\tau = 0$ then we obtain
\beq\lb{sosubreg0}
\mathcal{Y}_{1,0}=
 \{ v,w,z \mid z^{m-1} + v^2 + z w^2 = 0 \}.
\eeq
This is the defining equation of the 
Kleinian singularity of type $D_m$. As mentioned in the 
symplectic case, taking $\tau =0$ is the same as intersecting 
$\mathcal{S}_1$ with the nilpotent cone. So as before, we recover the 
result of Slodowy for this case.
Also, because of  $\lambda$-invariance, we know that  
the varieties $\mathcal{Y}_{1,\tau}$ are smooth
when $\tau \in R_{\ger{so}_{2m}(\Cc)}$. 
\end{example}

\subsection{Embedding into the Hilbert scheme}

As in Section \ref{sp} the defining equation
of $\mathcal{Y}_{n,\tau}$, modulo $\hat A$, is formally
the same as the defining equation of $\Gamma_{\tau}$. 
Moreover, if we are given polynomials $\hat V_i$ and $\hat W_i$
for $i=1,2$ as defined in the previous section, then both
$\hat V_1 - \hat V_2$ and $\hat W_1 - \hat W_2$ are polynomials
of degree no larger than $n-1$, while on the other hand the polynomial
$\hat A$ is  monic of degree $n$. Hence, the proof of 
Theorem \ref{sphilb} will carry over to the type-$D$ case 
and give us the following: 

\btm\lb{sohilb}
Let $2n+1 \leq m$ and 
$\tau \in R_{\ger{so}_{2m}(\Cc)} \subseteq \ger h / W$.
There is an open holomorphic embedding
\beq
j: \mathcal{Y}_{n,\tau} \to \mbox{Hilb}^n (\Gamma_{\tau}),
\eeq
where $\Gamma_{\tau}$
is the affine surface in $\Cc^3$ described by equation \eqref{sosubreg}.
\etm

Notice that $\Gamma_{\tau} \cong \Gamma_{c \tau}$ by the
map $(v,w,z) \mapsto (v,-w,z)$. Hence, we can 
consider $\tau$ to be in $\mbox{Sym}_m (\Cc)$ (via $\epsilon$) and
we will have a
surface $\Gamma_{\tau}$ that is well defined up to isomorphism.

\subsection{Type-$B$}

We have an obvious orbit preserving inclusion inclusion 
$i : \ger{so}_{2m}(\Cc) \into \ger{so}_{2m+1}(\Cc)$
given by
\[
i : Y \mapsto 
\left[
\begin{tabular}{c c}
0 & \\
  & Y
\end{tabular}
\right].
\]
Let  $X_n \in \ger{so}_{2m}(\Cc)$ be the nilpotent element defined
in \eqref{Xsoeven}, then $i(X_n)$ is a representative 
of the nilpotent orbit corresponding to the partition
$[2(m-n)-1,2n+1,1]$. Let us write $x_n = i(X_{n-1})$. 
Then $x_1$ is a subregular nilpotent element and in general
the codimension of the adjoint orbit of $x_n$ is $m +2n$.

The root system in the type-$B$ case is 
$\{ \pm e_i \pm e_j , \pm e_j \mid 1 \leq i,j \leq m, i \neq j \}$ with 
the positive roots being taken to be
$\{ e_i \pm e_j , e_k \mid 1 \leq i < j \leq m, 1 \leq k \leq m \}$.
The $\alpha$-root space is spanned by the vector $X_{\alpha}$ given by
\begin{align}
X_{e_i - e_j} &= E_{i+1,j+1} - E_{j+m+1, i+m+1}\\
X_{e_i + e_j} &= E_{i+1,j+m+1} - E_{j+1, i+m+1} \phantom{ajaj} (i<j)  \\
X_{-e_i - e_j} &= E_{i+m+1,j+1} - E_{j+m+1,i+1} \phantom{ajaj} (i<j)  \\
X_{e_i} &= E_{1,i+m+1} - E_{i+1,1}\\
X_{-e_i} &= E_{1,i+1} - E_{i+m+1,1}.
\end{align} 

We define the following affine subspace
of $\ger{so}_{2m+1}(\Cc)$:
\beq
\mathcal{S}_n^B = i(\mathcal{S}_{n-1}) + \mathcal{V}_n
\eeq
where $\mathcal{S}_{n-1}$ is the slice
of $\mathcal{O}_{X_n}$ defined in \eqref{soslice} and
$\mathcal{V}_n \subseteq \ger{so}_{2m+1}(\Cc)$ 
is the 2-dimensional vector subspace 
defined by
\beq
\mathcal{V}_n = \{ a_0  X_{- e_1} + d_0 X_{- e_{m-n+1}}  \mid x,y \in \Cc \}.
\eeq
We denote an element of $\mathcal{S}_n^B$ by
$S=S(a_i , y_i , z_j , d_k)$.

\btm
The affine subspace $\mathcal{S}_n^B$ is $\lambda$-invariant and
transverse to the adjoint orbit of $x_n$. 
Moreover, the characteristic
polynomial of an element $S = S (a_i , y_i , z_j , d_k )$ in
$\mathcal{S}_n^B$ is given by
\beq\lb{somainodd}
\chi_S (t) =  t a(t) d(t) 
              + (-1)^m t b(t) b(-t) 
\eeq
where
\begin{align}
a(t) &= t^2 A(t) + (-1)^n a_0^2 \\
d(t) &= D(t) + (-1)^{m-n} d_0^2 / t^2 \\ \notag
b(t) &=  a_0 d_0/t + B(t) \notag
\end{align}
where the polynomails $A$,$B$, and $D$ are those defined 
in the type-$D$ case by equtaions 
\eqref{soA}--\eqref{soB} (take n-1 instead of n).
\etm

The first thing to note is that the terms with $t$ in the denominator
will vanish if we simplify. However, we can also just multiply 
\eqref{somainodd} by $t$ so that we can preserve the form of the equation.
Also, there is a fiber-preserving $\Z / 2$-action
on $\mathcal{S}_n^B$ given by 
taking $(a_0, d_0)$ to $(-a_0 , -d_0 )$. 

Let us denote $\chi^{-1}(\tau) \cap \mathcal{S}_n^B$ by
$\mathcal{Y}_{n,\tau}^B$.
If we let define $U$ and $V$ as we did in the previous two cases and
then reduce all the polynomials we see that
for any $\tau \in \mbox{Sym}_m (\Cc)$ the variety 
$\mathcal{Y}_{n,\tau}^B$ is defined by the equation
\[
t P_{\tau} (t) + \hat U (t)^2 + t \hat V (t)^2 = \hat a (t) (t \hat d(t)).
\]
The polynomial $a(t)$ is monic of degree $n$ and the methods in the 
proof of 
Theorem \ref{sphilb} will carry over here to give us
\btm\lb{sohilbodd}
Let $2n-1 \leq m$ and 
$\tau \in R_{\ger{so}_{2m+1}(\Cc)} \subseteq \ger h / W$.
There is an algebraic morphism
\beq
j: \mathcal{Y}_{n,\tau}^B \to \mbox{Hilb}^n (\Sigma_{\{ 0 , \tau\} })
\eeq
which is not injective, but has fibers given by the $\Z / 2$-action
$(a_0,d_0) \mapsto (-a_0,-d_0)$.
\etm

\begin{example}
Let's look at the subregular case where we have
\begin{align}
a(t) &= t^2 - a_0^2 \\
d(t) &= D(t) + (-1)^{m-1} d_0^2 / t^2 \\
b(t) &= a_0 d_0 / t -z_1 + (-1)^m t^m.
\end{align}
So the defining equation of $\mathcal{Y}_{1,\tau}^B$ will be
\begin{align}
t^2 P_{\tau} (t^2) &= (t^2 - a_0^2)(t^2 D(t) + (-1)^{m-1} d_0^2)\\
    & \phantom{adada}   + (-1)^m (a_0 d_0 - z_1 t + (-1)^m t^{m+1})
                     ( - a_0 d_0 - z_1 t + t^{m+1})  \notag \\
    &= (t^2 - a_0^2)(t^2 D(t)) + (-1)^{m-1} d_0^2 t^2
        + (-1)^m z_1^2 t^2 \lb{sosubodd}\\
    & \phantom{adada}    - ((-1)^m + 1) z_1 t^{m+2}
        + ((-1)^m - 1) a_0 d_0 t^{m+1}
        + t^{2m+2}. \notag
\end{align}
Let us define
\beq
q(t) = t^{2m} -2 \delta (m-1) a_0 d_0 t^{m-1} - 2 \delta (m) z_1 t^m
\eeq
so that after dividing \eqref{sosubodd} by $t^2$ and rearranging slightly 
we get
\beq
P_{\tau} (t^2) + (-1)^{m} d_0^2 - (-1)^m z_1^2 - q(t)
= (t^2 - a_0^2)D(t)
\eeq
Now if we know that $a_0$ is a root of the left-hand side of this equation
we can recover $D$ uniquely. So we see that  $\mathcal{Y}_{1,\tau}^B$
will consist of the set of all $a_0, d_0, z_1$ such that 
\beq
P_{\tau} (a_0^2) + (-1)^{m} d_0^2 - (-1)^m z_1^2 - q(a_0) = 0.
\eeq
Now, depending on if $m$ is even or odd we get
\begin{align}
q(a_0) &= a_0^{2m} - 2 z_1 a_0^m &&m = \mbox{even} \\
q(a_0) &= a_0^{2m} - 2 d_0 a_0^m &&m = \mbox{odd}
\end{align}
which means that the defining equation of $\mathcal{Y}_{1,\tau}^B$ will be
\begin{align}
P_{\tau}(a_0^2)  + d_0^2 - (z_1 - a_0^m)^2 = 0 &&m = \mbox{even} \\
P_{\tau}(a_0^2)  - (d_0 - a_0^m)^2  + z_1^2 = 0 &&m = \mbox{odd}
\end{align}
Thus, after a change of coordinates we can write
\beq
\mathcal{Y}_{1,\tau}^B = \{ u,v,z \mid P_{\tau} +u^2 +v^2 =0 \}
\eeq
which is a deformation of the $A_{2m-1}$ singularity. 
This agrees with Proposition 
\ref{Slodowy}. 
\end{example}

Consider the above example together with Theorem \ref{sohilbodd}. We see that 
the type-$B$ case doesn't follow the pattern set by the previous cases.
In Theorem \ref{sohilbodd} we are mapping $\mathcal{Y}_{n,\tau}^B$ 
into the Hilbert scheme of points on a type-$D$ surface, but
this surface \emph{does not} correspond to the 
one of type-$A$  given by the subregular
orbit. However, if we were able to obtain a map into the Hilbert scheme over
the $A_{2m-1}$ surface using the methods of Theorem \ref{sphilb}, then we
would  be getting something substantially less
than what Theorem \ref{sohilbodd} gives us. 
Namely, instead of modding out by the degree $n$
reduced polynomial $\hat a (t)$ we would have to mod out by the 
unreduced polynomial $a (t)$ which has degree $2n$. So we would
be embedding the $2n$-dimensional 
$\mathcal{Y}_{n,\tau}^B$ into the $4n$-dimensional 
$\mbox{Hilb}^{2n}(S_{\tau})$.


\end{document}